\documentclass[12pt,letterpaper]{article}
\usepackage[T1]{fontenc}
\usepackage[english]{babel}	
\usepackage{german}

\usepackage{graphicx}
\usepackage[a4paper]{geometry}
\usepackage{fullpage}
\usepackage{hyperref}

\usepackage{algorithm}
\usepackage{algpseudocode} 
\usepackage{amsmath}
\usepackage{amsthm}
\usepackage{amssymb}
\usepackage{mathpazo}
\usepackage{array}
\usepackage{nicefrac}

\theoremstyle{definition}
\newtheorem{defi}{Definition} 

\theoremstyle{plain}

\numberwithin{equation}{section}

\newcommand{\mN}{\mathbb{N}}
\newcommand{\mZ}{\mathbb{Z}}
\newcommand{\mP}{\mathbb{P}}

\newcommand{\mdiv}{\!\mid\!}

\newcommand{\copr}{\!\perp\!}

\newcommand{\vl}{\vrule width 2pt}

\begin{document}

\title{A short note on the computation\linebreak of the generalised Jacobsthal function\linebreak for paired progressions}
\author{Mario Ziller and John F. Morack}
\date{}

\maketitle

\begin{abstract}

Jacobsthal's function was recently generalised for the case of paired\linebreak progressions. It was proven that a specific bound of this function is sufficient for the truth of Goldbach's conjecture and of the prime pairs conjecture as well. We extended and  adapted algorithms described for the computation of the common\linebreak Jacobsthal function, and computed respective function values of the paired\linebreak Jacobsthal function for primorial numbers for primes up to 73. All these values fulfil the conjectured specific bound. In addition to this note, we provide a detailed review of the algorithmic approaches and the complete computational results in ancillary files.\\\\
\end{abstract}

\section{Introduction}

Henceforth, we denote the set of integral numbers by $\mZ$ and the set of natural\linebreak numbers, i.e. positive integers, by $\mN$. Let $\mP=\{p_i\mid i\in\mN\}$ be the set of prime numbers with $p_1=2$.\\

The commonly known Jacobsthal function $j(n)$ is defined to be the smallest natural number $m$, such that every sequence of $m$ consecutive integers contains at least one integer coprime to $n$ \cite{Jacobsthal_1960_I, Erdos_1962}.

\begin{defi} \label{Jacobsthal}
For $n\in\mN$, the Jacobsthal function $j(n)$ is defined as
$$j(n)=\min\ \{m\in\mN\mid\forall\ a\in\mZ\ \exists\ q\in\{1,\dots,m\}:a+q\copr n\}.$$
\end{defi}

Jacobshal remarked that the entire function is determined by its values for\linebreak products of distinct primes \cite{Jacobsthal_1960_I}. The particular case of primorial numbers is therefore\linebreak of great interest. The function values at these points characterise the maximum growth of the function.\\

\pagebreak

The Jacobsthal function was recently generalised for the case of paired progressions. A paired progression \cite{Ziller_Morack_2017} is defined as an ordered sequence of consecutive pairs of integers of a definite length.

\begin{defi}
Let $a,b\in\mZ$ and $m\in\mN$. A paired progression $\langle a,b\rangle_m$ is defined as
$$\langle a,b\rangle_m=\{(a+i,b+i)\in\mZ^2\mid i=1,\dots,m\}.$$
\end{defi}

The generalisation $j_2(n)$ of the Jacobsthal function with respect to paired\linebreak progressions is called paired Jacobsthal function. It is defined to be the smallest natural\linebreak number $m$, such that every paired progression of length $m$, where the difference of its pair elements is even, contains at least one pair coprime to $n$ \cite{Ziller_Morack_2017}.

\begin{defi}
For $n\in\mN$, the paired Jacobsthal function $j_2(n)$ is defined as
\begin{equation*} \begin{split}
j_2(n)=\min\ \{m\in\mN\mid\ &\forall\ (a,b)\in\mZ^2\ with\ 2\mdiv (b-a):\\
&\exists\ q\in\{1,\dots,m\}:n\copr (a+q,b+q)\}.
\end{split} \end{equation*}
\end{defi}

This function is also determined by its values for products of distinct primes \cite{Ziller_Morack_2017}, and the case of primorial numbers is most interesting. The function values at these points generate a separate function.

\begin{defi}
For $n\in\mN$, the  primorial paired Jacobsthal function $h_2(n)$ is defined as
\begin{equation*} \begin{split}
h_2(n)=\min\ \{m\in\mN\mid\ &\forall\ (a,b)\in\mZ^2\ with\ 2\mdiv (b-a):\\
&\exists\ q\in\{1,\dots,m\}\ \forall\ i\in\{1,\dots,n\}:p_i\copr (a+q,b+q)\}.
\end{split} \end{equation*}
\end{defi}

It was proven that the conjectured upper bound $h_2(n)<p_n^2-p_n,\ n\ge 3$ of this function represents a sufficient condition for the truth of Goldbach's conjecture and of the prime pairs conjecture as well \cite{Ziller_Morack_2017}.\\\\

\section{Computation and results}

Several algorithmic concepts were successfully applied to the computation of\linebreak Jacobsthal's common function for primorials \cite{Ziller_Morack_2016}. For the computation of the primorial paired Jacobsthal function, we extended and adapted these approaches for the case of paired progressions. A detailed analysis of how to compute $h_2(n)$ is provided in an ancillary file.\\

We computed the entire values of the function $h_2(n)$ for $n\le21$ or $p_n\le73$, respectively.  The results for $3\le n\le21$ fulfil the conjectured bound  $h_2(n)<p_n^2-p_n$ \cite{Ziller_Morack_2017}. The data are listed in table 1, including the maximum processed prime $p_n$, the values of $h_2(n)$, and the conjectured bound $p_n^2-p_n$.\\

\pagebreak
\newpage

\vspace*{1mm}

\begin{table}[!h]
  \centering
\begin{tabular}{cc}	
\begin{tabular}{!\vl rrr|r !\vl}
  \noalign{\hrule height 2pt} 
\rule{0pt}{14pt}$n$&$p_n$&$h_2(n)$&$p_n^2$-$p_n$\\[2pt]
  \noalign{\hrule height 2pt} 
  \rule{0pt}{14pt}
1&2&2&2\\
2&3&6&6\\
3&5&18&20\\
4&7&30&42\\
5&11&66&110\\
6&13&150&156\\
7&17&192&272\\
8&19&258&342\\
9&23&366&506\\
10&29&450&812\\
11&31&570&930\\
  [2pt]
  \noalign{\hrule height 2pt}
\end{tabular}
&	 \hspace{5pt}
\begin{tabular}{!\vl rrr|r !\vl}
  \noalign{\hrule height 2pt} 
\rule{0pt}{14pt}$n$&$p_n$&$h_2(n)$&$p_n^2$-$p_n$\\[2pt]
  \noalign{\hrule height 2pt} 
  \rule{0pt}{14pt}
12&37&708&1332\\
13&41&894&1640\\
14&43&1044&1806\\
15&47&1284&2162\\
16&53&1422&2756\\
17&59&1656&3422\\
18&61&1902&3660\\
19&67&2190&4422\\
20&71&2460&4970\\
21&73&2622&5256\\
&&&\\
  [2pt]
  \noalign{\hrule height 2pt}
\end{tabular}
\end{tabular}\\
 \caption{Results of computation.}
\end{table}


\bigskip

In addition to this note, we provide several ancillary files. The file \dq full\_details.pdf\dq\ includes a detailed description of theory, algorithms, and further results. This paper\linebreak widely follows the structure of our former paper \cite{Ziller_Morack_2016}. Furthermore, it outlines the\linebreak content of all other anxillary data files representing exhaustive lists of the computational results.\\

\subsubsection*{Contact}
marioziller@arcor.de\\
axelmorack@live.com\\

\bibliography{References}     

\end{document}